\documentclass[12pt,letterpaper]{article}
\usepackage{amsfonts}
\usepackage{amssymb}
\usepackage{latexsym}  % latex2e symbol package
\usepackage{pifont}  % so we can use ZapfDingbats

\let\SavedRightarrow=\Rightarrow
\usepackage{marvosym}
\let\Rightarrow=\SavedRightarrow

\usepackage{amsmath}   % so we can use \tag , etc.

\newcommand\PP{{\mathcal P}}
\newcommand\BB{{\mathcal B}}
\newcommand\FF{{\mathcal F}}
\newcommand\GG{{\mathcal G}}

\newcommand\cchi{{\raise 2 pt \hbox{$\chi$}}}
\newcommand{\MS}{\mathsf{MS}}  % Menger Sponge

\newcommand\ro{\mathsf{ro}}  % regular open sets
\newcommand\acl{\mathsf{acl}}  % almost clopen

\newcommand\cl{\mathrm{cl}}   % closure
\newcommand\intr{\mathrm{int}}   % interior
\newcommand\lh{\mathrm{lh}}   % length

\newcommand\cat{^{\mathord{\frown}}}  % for string concatenation

\newcommand\compl{^{\mbox{\tiny $\complement$}}} % boolean complement

\newcommand\rest{\upharpoonright}     % as a relation
\newcommand\res{\mathord {\upharpoonright}}  % less space around it

\newcommand\ind{\mathrm{ind}}  % inductive dimension
\newcommand\trind{\mathrm{trind}}  % transfinite inductive dimension
\newcommand\iv{^{-1}} % inverse

\newcommand\onto{\twoheadrightarrow}

\newcommand\eop{\raisebox{-2pt}{\mbox{\Huge \Smiley}}}

{\end{enumerate}}
\newenvironment{itemizz}{\begin{itemize}\setlength{\itemsep}{-1mm}} %
{\end{itemize}}                              

\newtheorem{theorem}{Theorem}[section]
\newtheorem{definition}[theorem]{Definition}
\newtheorem{lemma}[theorem]{Lemma}
\newtheorem{remark}[theorem]{Remark}
\newtheorem{conditions}[theorem]{Conditions}

\newenvironment{proof}{{\bf Proof.}}{\eop\medskip}
\newenvironment{proofof}[1]{\medskip \textbf{Proof of #1.}}{\eop\medskip}

\begin{document}

\title{One Dimensional Locally Connected S-spaces
}

\author{Joan E. Hart\footnote{University of Wisconsin, Oshkosh,
WI 54901, U.S.A.,
\ \ hartj@uwosh.edu}
\  and
Kenneth Kunen\footnote{University of Wisconsin,  Madison, WI  53706, U.S.A.,
\ \ kunen@math.wisc.edu}
\thanks{Both authors partially supported by NSF Grant 
DMS-0456653.}
}

\maketitle

\begin{abstract}
We construct, assuming Jensen's principle $\diamondsuit$, a  one-dimensional
locally connected hereditarily separable
continuum without convergent sequences.
\end{abstract}

\section{Introduction} 
\label{sec-intro}
All topologies discussed in this paper are assumed to be Hausdorff.
A \emph{continuum} is any compact connected space.
A \emph{nontrivial convergent sequence} is a convergent
$\omega$--sequence of distinct points.
As usual, $\dim(X)$ is the covering dimension of $X$; for details, see
Engelking \cite{ENGd}.   ``HS'' abbreviates ``hereditarily separable''.
We shall prove:

\begin{theorem}
\label{thm-main}
Assuming $\diamondsuit$, there is a locally connected
HS continuum $Z$ such that
$\dim(Z) = 1$ and $Z$ has no nontrivial convergent sequences.
\end{theorem}

Note that points in $Z$ must have uncountable character,
so that $Z$ is not hereditarily Lindel\"of; thus, $Z$ is an S-space.

Spaces with some of these features are well-known from the literature.
A compact F-space has no nontrivial convergent sequences.
Such a space can be a continuum; for example, the
\v Cech remainder $\beta[0,1)\setminus[0,1)$ is connected,
although not locally connected;
more generally, no infinite compact F-space can be either
locally connected or HS.
In \cite{vM}, van Mill constructs, under the Continuum Hypothesis,
a locally connected continuum
with no nontrivial convergent sequences.
Van Mill's example, constructed as an inverse limit
of Hilbert cubes, is infinite dimensional.
Here, we shall replace the Hilbert cubes by one-dimensional
Peano continua (i.e., connected, locally connected, compact
metric spaces) to obtain a one-dimensional limit space.
Our $Z = Z_{\omega_1}$ will be the limit
of an inverse system $\langle Z_\alpha : \alpha < \omega_1\rangle$.
Each $Z_\alpha$ will be a copy of the
\emph{Menger sponge} \cite{Me} (or Menger curve) $\MS$; this
one-dimensional Peano continuum has homogeneity properties
similar to those of the Hilbert cube.
The basic properties of $\MS$ are summarized in
Section \ref{sec-sponge}, and
Theorem \ref{thm-main} is proved in Section~\ref{sec-inv}.

In \cite{vM}, as well as in earlier work by
Fedorchuk \cite{FED} and van Douwen and Fleissner \cite{vDF},
one kills all possible
nontrivial convergent sequences in  $\omega_1$ steps.
Here, we focus primarily on obtaining an S-space,
modifying the construction of the original Fedorchuk S-space \cite{FED-S};
we follow the exposition in \cite{DK}, where the lack of convergent
sequences occurs only as an afterthought.

We do not know whether one can obtain $Z$  so that it satisfies
Theorem \ref{thm-main} with
the stronger property $\ind(Z) = 1$;  that is, 
the open $U \subseteq Z$ with $\partial U$ zero-dimensional form a base.
In fact, we can easily modify our construction to ensure that
$1=\dim(Z)< \ind(Z) = \infty$; this will hold because
(as in \cite{DK}) we can give $Z$ the additional property that
all perfect subsets are $G_\delta$ sets;
see  Section \ref{sec-rem} for details.

We can show that a $Z$ satisfying Theorem \ref{thm-main} cannot
have the property that the open $U \subseteq Z$ with $\partial U$
scattered form a base;
see  Theorem \ref{thm-acl-not-base} in Section \ref{sec-acl}.
This strengthening of  $\ind(Z) = 1$ is satisfied by some
well-known Peano continua.  It is also satisfied by the space
produced in \cite{HK2} under  $\diamondsuit$
by an inductive construction related to the one we describe here,
but the space of \cite{HK2}
was not locally connected, and it had nontrivial
convergent sequences (in fact, it was hereditarily Lindel\"of).

\section{On Sponges}
\label{sec-sponge}
The Menger sponge $\MS$ \cite{Me} is obtained by drilling holes
through the cube $[0,1]^3$, analogously to the way that
one obtains the middle-third Cantor set
by removing intervals from $[0,1]$.
The paper of  Mayer,  Oversteegen, and Tymchatyn \cite{MOT} has
a precise definition of $\MS$ and
discusses its basic properties.

In proving theorems about $\MS$, one often refers not to its
definition, but to the following theorem of R. D. Anderson \cite{AND2,AND1}
(or, see \cite{MOT}), which  characterizes $\MS$.
This theorem will be used to verify inductively that $Z_\alpha \cong \MS$.
The fact that $\MS$ satisfies the stated conditions is easily
seen from its definition, but it is not trivial
to prove that they characterize $\MS$.

\begin{theorem}
\label{thm-MS-char}
$\MS$ is, up to homeomorphism, the only one-dimensional Peano continuum with no
locally separating points and no non-empty planar open sets.
\end{theorem}

Here, $C \subseteq X$ is  \emph{locally separating} iff,
for some connected open $U \subseteq X$,
the set $U\setminus C$ is not connected.
A point $x$ is locally separating iff $\{x\}$ is.
This notion is applied in
the Homeomorphism Extension Theorem
of  Mayer,  Oversteegen, and Tymchatyn \cite{MOT}:

\begin{theorem}
\label{thm-hom-ext}
Let $K$ and $L$ be closed, non-locally-separating subsets of
$\MS$ and let $h: K \onto L$ be a homeomorphism.
Then $h$ extends to a homeomorphism of $\MS$ onto itself.
\end{theorem}

The non-locally-separating sets have the following
closure property of Kline \cite{Kl} (or, see Theorem 2.2 of \cite{MOT}):

\begin{theorem}
\label{thm-seq-sep}
Let $X$ be compact and locally connected, and let
$K = \bigcup\{K_i : i \in \omega\}$,
where $K$ and the $K_i$ are closed subsets of $X$.
If $K$ is locally separating then some $K_i$ is locally separating.
\end{theorem}

For example, these results imply that in $\MS$, all convergent sequences are equivalent.
More precisely,
points in $\MS$ are not locally separating, so if
$\langle x_i : i \in \omega\rangle$ converges to $x_\omega$,
then $\{x_i : i \le \omega\}$ is not locally separating.
Thus, if $\langle s_i \rangle$ and $\langle t_i\rangle$ 
are nontrivial convergent sequences in $\MS$,
with limit points $s_\omega$ and $t_\omega$, respectively,
then there is a homeomorphism of $\MS$ onto itself that
maps $s_i$ to $t_i$ for each $i \le \omega$.

The following  consequence of Theorem \ref{thm-MS-char} was noted by
Prajs \cite{Pr} (see p.~657).

\begin{lemma}
\label{lemma-collapse-arc}
Let $J \subseteq \MS$ be a non-locally-separating
arc and obtain $\MS/ J$ by collapsing
$J$ to a point.  Then $\MS/ J \cong \MS$ and the natural map
$\pi : \MS \onto \MS/J$ is monotone.
\end{lemma}

Here,
a map $f : Y \onto X$ is called \emph{monotone} iff each
$f\iv\{x\}$ is connected; so, the monotonicity in Lemma \ref{lemma-collapse-arc}
is obvious.
When $X,Y$ are compact, monotonicity implies that 
$f\iv(U)$ is connected whenever $U$ is a connected open or closed
subset of $X$.

We shall use these results to
show that the property of being a Menger sponge will
be preserved at the limit stages of our construction:

\begin{lemma}
\label{lemma-limit}
Suppose that $\gamma$ is a countable limit ordinal and $Z_\gamma$
is an inverse limit of $\langle Z_\alpha : \alpha < \gamma \rangle$,
where all bonding maps $\pi^\beta_\alpha$ are monotone and each $Z_\alpha \cong \MS$.
Then $Z_\gamma \cong \MS$.
\end{lemma}
\begin{proof}
We  verify the conditions of Theorem \ref{thm-MS-char}.
$\dim(Z_\gamma) = 1$, since this property is preserved by
inverse limits of compacta,
and $Z_\gamma$ is locally connected because the
$\pi^\beta_\alpha$ are monotone.  So,  we need to verify that
$Z_\gamma$ has no
locally separating points and no non-empty planar open sets.

Suppose that $q \in Z_\gamma$ is locally separating;
so we have a connected neighborhood $U$ of $q$ with $U \backslash \{q\}$
not connected.
Shrinking $U$, we may assume that $U = (\pi^\gamma_\alpha)\iv(V)$, where
$\alpha < \gamma$ and $V$ is open and connected in $Z_\alpha$.
Since $Z_\alpha \cong \MS$, $\pi^\gamma_\alpha(q)$ is not locally separating,
so $V \backslash \{ \pi^\gamma_\alpha(q) \}$ is connected.
Then, since $\pi^\gamma_\alpha$ is monotone,
$(\pi^\gamma_\alpha)\iv ( V \backslash \{ \pi^\gamma_\alpha(q)\} ) =
U \setminus (\pi^\gamma_\alpha)\iv  \{ \pi^\gamma_\alpha(q)\}$
is connected.  The same argument shows that
$ U \setminus (\pi^\gamma_\beta)\iv  \{ \pi^\gamma_\beta(q)\}$ is
connected whenever $\alpha \le \beta < \gamma$.
But then $U \backslash \{q\} = \bigcup \{
U \setminus (\pi^\gamma_\beta)\iv  \{ \pi^\gamma_\beta(q)\} :
\alpha \le \beta < \gamma \}$ is connected also.

Suppose that $U \subseteq Z_\gamma$ is open and non-empty;
we show that $U$ is not planar.  Shrinking $U$,
we may assume that $U = (\pi^\gamma_\alpha)\iv(V)$, where
$\alpha < \gamma$ and $V$ is open in $Z_\alpha$.
Since $Z_\alpha \cong \MS$, there is a $K_5$ set $F \subseteq V$;
that is, $F$ consists of $5$ distinct points $p_0, p_1, p_2, p_3, p_4$
together with arcs $J_{i,j}$ with endpoints $p_i,p_j$ for $0 \le i < j < 5$,
where the sets $J_{i,j} \setminus \{p_i, p_j\}$, for $0 \le i < j < 5$,
are pairwise disjoint.  Now $F$ is not planar, and,
one can show that $(\pi^\gamma_\alpha)\iv (F)$ is not planar either.
To do this, use the fact that $\pi^\gamma_\alpha$ is monotone, so that
the sets $(\pi^\gamma_\alpha)\iv \{p_i\}$ and
$(\pi^\gamma_\alpha)\iv (J_{i,j})$ are all continua.
\end{proof}

The following terminology was used also in
the exposition in \cite{DK} of the Fedorchuk S-space:

\begin{definition}
\label{def-strong-limit}
Let $\FF$ be a family of subsets of $X$.
Then $x \in X$ is a \emph{strong limit point} of $\FF$ iff
for all neighborhoods $U$ of $x$, there is an $F \in \FF$
such that $F \subseteq U$ and $x \notin F$.
\end{definition}

In practice, we shall only use this notion when the elements of $\FF$
are closed.  If all elements of $\FF$
are singletons, this reduces to the usual notion of a point
being a limit point of a set of points.

The map $\pi^{\alpha+1}_\alpha : Z_{\alpha+1} \onto Z_\alpha$
will always be obtained by collapsing a non-locally-separating
arc in $Z_{\alpha+1}$ to a point.   We obtain it using:

\begin{lemma}
\label{lemma-succ-stage}
Assume that $X \cong \MS$ and that for $n \in \omega$,
$\FF_n$ is a family of non-locally-separating closed subsets
of $X$.  Fix $t \in X$ such that $t$ is a strong limit
point of each $\FF_n$.
Then there is a $Y \cong \MS$ and a monotone $\pi : Y \onto X$
such that $\pi\iv\{t\}$ is a non-locally-separating arc in $Y$,
$|\pi\iv\{x\}| = 1$ for all $x \ne t$,
and, for each $n$ and each $y \in \pi\iv\{t\}$:
$y$ is a strong limit point of $\{\pi\iv(F) : F \in \FF_n\}$.
\end{lemma}
\begin{proof}
First, let $\{A_n : n \in \omega\}$
partition $\omega$ into disjoint infinite sets.
In $X$, choose disjoint closed $F_i \not\ni t$ for $i \in \omega$ such that
$F_i \in \FF_n$ whenever $i \in A_n$, and such that
every neighborhood of $t$ contains all but finitely many of the $F_i$.
Let $L = \{t\} \cup \bigcup_i F_i$.  Then $L$ is closed
and non-locally-separating by Theorem \ref{thm-seq-sep}.

Now, in $\MS$, let $J$ be any non-locally-separating arc.
Choose disjoint closed non-locally separating sets $G_i$ for $i \in \omega$
such that each $G_i \cong F_i$, every neighborhood of $J$
contains all but finitely many $G_i$, each $G_i \cap J = \emptyset$,
and for each $n$ and each $y \in J$:
$y$ is a strong limit point of $\{G_i : i \in A_n\}$.

Let $\sigma : \MS \onto \MS/J$ be the usual projection,
and let $[J]$ denote the point to which $\sigma$ collapses the set $J$.
Then $\MS/ J \cong \MS$ by Lemma \ref{lemma-collapse-arc}.
In $\MS/ J$, let $K = \{[J]\} \cup \bigcup\{\sigma(G_i) : i \in \omega\}$.
Let $h: K \onto L$ be a homeomorphism such that $h([J]) = t$
and each $h( \sigma(G_i) ) = F_i$.  By Theorem \ref{thm-hom-ext},
$h$ extends to a homeomorphism $\widetilde h : \MS/J \onto X$.

Now, let $Y = \MS$ and let $\pi = \widetilde h \circ \sigma$.
\end{proof}

The next lemma will simplify somewhat the description of our inverse limit:

\begin{lemma}
\label{lemma-get-in-cube}
In Lemma \ref{lemma-succ-stage}, we may obtain $Y \subseteq X \times [0,1]$,
with $\pi : Y \onto X$ the natural projection.
\end{lemma}
\begin{proof}
Start with any $Y,\pi,t$ satisfying Lemma \ref{lemma-succ-stage},
and let $J := \pi^{-1}\{t\}$.
Apply the Tietze Extension Theorem to fix $f : Y \onto [0,1]$ such 
that $f \res J : J \onto [0,1]$ is a homeomorphism.
Then $y \mapsto (\pi(y), f(y))$ is one-to-one on $Y$,
and hence
$\widetilde Y := \{(\pi(y), f(y)) : y \in Y \} \subseteq X \times [0,1]$
satisfies Lemma \ref{lemma-get-in-cube}.
\end{proof}

The following additional property of our $\pi$ will be useful:

\begin{lemma}
\label{lemma-preserve-non-loc-sep}
Let $t$ and $\pi : Y \onto X$ be as in Lemma
\ref{lemma-succ-stage} or \ref{lemma-get-in-cube}.
Assume that $H \subseteq X$ is closed and nowhere dense
and not locally separating.  Then
$\pi\iv(H) \subseteq Y$ is closed and nowhere dense
and not locally separating.
\end{lemma}
\begin{proof}
$\pi\iv(H)$ is closed and nowhere dense because $\pi$ is continuous
and irreducible.
Also note that $\pi\iv(H)$ is not locally separating 
if either $H = \{t\}$ (trivially) or $t \notin H$
(because $\pi$ is a homeomorphism in a neighborhood of  $\pi\iv(H)$).

Next, note that every closed  $K \subseteq H$ 
is non-locally-separating in $X$:  If not, let $U \subseteq X$
be connected and open with $U \backslash K$ not connected, so that
$U \backslash K = W_0 \cup W_1$, where the $W_i$ are
open in $X$, non-empty, and disjoint.  Then
$U \backslash H = W_0\backslash H \,\cup\, W_1\backslash H$,
but $H$ is not locally separating, so one of the
$W_i\backslash H = \emptyset$, so 
$W_i \subseteq H$, contradicting $H$ being nowhere dense.

Now, let $H = \bigcup_{n \in \omega} K_n$, where each $K_n$ is closed
and  either  $K_n = \{t\}$  or $t \notin K_n$.
Then $\pi\iv(H) = \bigcup_n \pi\iv(K_n)$, which is  not locally separating
by Theorem \ref{thm-seq-sep}.
\end{proof}

\section{The Inverse Limit}
\label{sec-inv}
We shall obtain our space $Z = Z_{\omega_1}$ as an inverse limit
of a sequence $\langle Z_\alpha : \alpha < \omega_1\rangle$.
As with many such constructions, it is somewhat simpler to view the
$Z_\alpha$ concretely as subsets of cubes, so that the bonding maps
are just projections.  Thus, we shall have:

\begin{conditions}
We obtain $Z_\alpha$ for $\alpha \le \omega_1$ and
$\pi^\beta_\alpha, \sigma^\beta_\alpha$
for $\alpha \le \beta \le \omega_1$ such that:
\begin{itemizz}
\item[C1.] Each $Z_\alpha$ is a closed subset of $\MS \times [0,1]^\alpha$,
and $Z_0 = \MS$.
\item[C2.] For $\alpha \le \beta \le \omega_1$,
$\pi^\beta_\alpha : \MS \times [0,1]^\beta \onto \MS \times [0,1]^\alpha$
is the natural projection.
\item[C3.] $\pi^\beta_\alpha(Z_\beta) = Z_\alpha$ whenever 
$\alpha \le \beta \le \omega_1$.
\item[C4.] $Z_\alpha$ is homeomorphic to $\MS$ whenever $\alpha < \omega_1$.
\item[C5.] The maps $\sigma^\beta_\alpha  :=
\pi^\beta_\alpha \res Z_\beta : Z_\beta \onto Z_\alpha$,
for $\alpha \le \beta \le \omega_1$, are monotone.
\end{itemizz}
\end{conditions}

Using (C1,C2,C3), the construction is determined
at limit ordinals; (C4) is preserved by Lemma \ref{lemma-limit} and (C5).
It remains to explain how,
given $Z_\alpha$ for $\alpha < \omega_1$, we obtain
$Z_{\alpha + 1} \subseteq Z_\alpha \times [0,1]$;
as usual, we identify $\MS \times [0,1]^{\alpha + 1}$ with
$\MS \times [0,1]^{\alpha} \times [0,1]$.

We now add:

\begin{conditions}
We have $q_\alpha^\xi$ and $t_\alpha$ for $\xi < \alpha < \omega_1$ such that:
\begin{itemizz}
\item[C6.]  Each $\langle q_\alpha^\xi : \xi < \alpha \rangle$
is a sequence of points in  $\MS \times [0,1]^\alpha$.
\item[C7.]  Whenever 
 $\langle q^\xi : \xi < \omega_1 \rangle$
is any sequence of points in  $\MS \times [0,1]^{\omega_1}$,
$\{ \alpha < \omega_1 : \forall \xi < \alpha
\, [\pi^{\omega_1}_\alpha(q^\xi) = q_\alpha^\xi ]\} $
is stationary.
\item[C8.] Whenever $\alpha < \beta  \le \omega_1$ and $z \in Z_\alpha$:
If $q_\alpha^\xi \in Z_\alpha$ for all $\xi < \alpha$ and $z$ is a limit point
of $\langle q_\alpha^\xi : \xi < \alpha \rangle$, then
all points of $(\sigma^\beta_\alpha)\iv \{z\}$ are strong
limit points of 
$\{ (\sigma^\beta_\alpha)\iv \{ q_\alpha^\xi \} :
\xi < \alpha \}$.
\item[C9.]  $t_\alpha \in Z_\alpha$, and for all $z \in Z_\alpha$:
$(\sigma^{\alpha+1}_\alpha)\iv\{z\}$ is a singleton
if $z \ne t_\alpha$ and a non-locally-separating arc if $z = t_\alpha$.
\item[C10.]  $t_\alpha = q^0_\alpha$ whenever $\alpha > 0$
and $q^0_\alpha \in Z_\alpha$.
\end{itemizz}
\end{conditions}

\begin{proofof}{Theorem \ref{thm-main}}
The fact that one may obtain
(C1 -- C10) has already been outlined above.
(C6,C7) are possible by $\diamondsuit$, and (C10) is just a definition.
(C8,C9) are obtained by induction on $\beta$.
For the successor step, we must obtain $Z_{\beta+1}$ from
$Z_\beta$ using Lemmas \ref{lemma-succ-stage} and \ref{lemma-get-in-cube}.
Here, $X = Z_\beta$, $Y = Z_{\beta+1}$, and $t = t_\beta$; the 
$\FF_n$ list all sets of the form
$\FF^\beta_\alpha := \{  (\sigma^\beta_\alpha)\iv \{ q_\alpha^\xi \} :
\xi < \alpha \ \&\ q_\alpha^\xi \in Z_\alpha \}$ such that $\alpha \le \beta$
and $t_\beta$ is a strong limit point of $\FF^\beta_\alpha$.
Observe that (C8) for $(\alpha, \beta+1)$ is immediate from
(C8) for $(\alpha, \beta)$ \emph{except} for the points of $Z_{\beta+1}$
in $(\sigma^{\beta+1}_\beta)\iv \{t_\beta\}$.
Also observe that in order to apply 
Lemmas \ref{lemma-succ-stage} and \ref{lemma-get-in-cube}, we must check
by induction on $\beta$, using Lemma \ref{lemma-preserve-non-loc-sep},
that the sets $ (\sigma^\beta_\alpha)\iv \{ q_\alpha^\xi \}$ are
non-locally-separating (and nowhere dense) in $Z_\beta$.

Note that $\cchi(z,Z) = \aleph_1$ for all $z \in Z$; this follows from
(C9,C10) and the fact, using (C7), that
$\{ \alpha < \omega_1 : \pi^{\omega_1}_\alpha(z) = t_\alpha \} $
is unbounded in $\omega_1$.

$Z$ is HS by (C6,C7,C8,C1,C2,C3):  If not, suppose that 
$\langle q^\xi : \xi < \omega_1 \rangle$ is left-separated in $Z$.
As in \cite{DK}, we get a club $C \subset \omega_1$ such that
for all $\alpha\in C$, 
\[
\cl \left\{ \sigma^{\omega_1}_\alpha(q^\xi) : \xi < \alpha \right\} \supseteq
 \left\{ \sigma^{\omega_1}_\alpha(q^\xi) : \xi < \omega_1\right\} \ \ .
\]
Fix $\alpha\in C$ such that
$ \forall \xi < \alpha \, [\sigma^{\omega_1}_\alpha(q^\xi) = q_\alpha^\xi ] $.
Let $z = \sigma^{\omega_1}_\alpha(q^\alpha)$.
Applying (C8) with $\beta= \omega_1$, we have in $Z$:
all points of $(\sigma^{\omega_1}_\alpha)\iv \{z\}$ are strong
limit points of 
$\{(\sigma^{\omega_1}_\alpha)\iv \{ q_\alpha^\xi \} :
\xi < \alpha \}$.  In particular, $q^\alpha$ is a limit point
of $\langle q^\xi : \xi < \alpha \rangle$, contradicting ``left-separated''.

Similarly, $Z$ has no non-trivial convergent 
sequences:  Suppose that $q^n \to q^\omega$ in $Z$,
where the $q^\xi$ for $\xi \le \omega$ are distinct.
Let $q^\xi =  q^\omega$ when $\omega < \xi < \omega_1$,
and apply (C7) to get $\alpha$ with $\omega < \alpha < \omega_1$ such that the
$\sigma^{\omega_1}_\alpha(q^\xi)$ for $\xi \le \omega$ are distinct points and
$ \forall \xi < \alpha \, [\sigma^{\omega_1}_\alpha(q^\xi) = q_\alpha^\xi ]$.
Let $z = \sigma^{\omega_1}_\alpha(q^\omega)$.
Then all points of $(\sigma^{\omega_1}_\alpha)\iv \{z\}$ are strong
limit points of
$\{ (\sigma^{\omega_1}_\alpha)\iv \{ q_\alpha^\xi \} : \xi < \alpha \}$
and hence also of
$\{ (\sigma^{\omega_1}_\alpha)\iv \{ q_\alpha^n \} : n < \omega \}$.
So, all points of  $(\sigma^{\omega_1}_\alpha)\iv \{z\}$ are
limit points of $\{q^n : n \in \omega\}$.  Since
$\{q^\omega\} \subsetneqq  (\sigma^{\omega_1}_\alpha)\iv \{z\}$
(by $\cchi(q^\omega,Z) = \aleph_1$), we contradict $q^n \to q^\omega$.
\end{proofof}

\section{The Almost Clopen Algebra}
\label{sec-acl}

We show here (Theorem \ref{thm-acl-not-base})
that a space $Z$ satisfying Theorem \ref{thm-main} cannot have
a base of open sets with scattered boundaries; equivalently
(because there are no nontrivial convergent sequences) with finite boundaries.
We first note that if there were such a base,
we could take the basic open sets $U$ to be regular,
since $\partial(\intr(\cl(U))) \subseteq \partial U$.
To simplify notation, we define:

\begin{definition}
$\ro(X)$ denotes the algebra of regular open subsets of $X$,
and $\acl(X)$
\textup(the \emph{almost clopen} sets\textup)
denotes the family of regular open sets $U$
such that $\partial U$ is finite.
For $U \in \ro(X)$, let
$U\compl$ denote the boolean complement $(X \backslash U)^\circ$.
\end{definition}

Note that $\acl(X)$ is a boolean subalgebra of $\ro(X)$:
If $U \in \acl(X)$ and $W = U\compl$, then $\partial W = \partial U$,
so $W \in \acl(X)$.
Also, if $U,V \in \acl(X)$ and $W = U \wedge V = U \cap V \in \ro(X)$,
then $W \in \acl(X)$ because
$\partial(W) \subseteq \partial(U) \cup \partial(V)$.

In a locally connected space,
the connected components of an open set $U$ are open;
if $V$ is any such component, then $\partial V \subseteq \partial U$
(because $V$ is relatively clopen in $U$), so $V \in \acl(X)$
whenever $U \in \acl(X)$.  Thus,

\begin{lemma}
If $X$ is locally connected and $\acl(X)$ is a local base
at $p \in X$, then $\{ U \in \acl(X) : p \in U \ \&\  U$ is connected$ \, \}$
is also a local base at $p$.
\end{lemma}

Various LOTS sums have bases of almost clopen sets.
This is true, for example, for any
compact hedgehog consisting of a central point plus
arbitrarily many LOTS spines.  The assumption of no
convergent sequences, however,
puts some restrictions on the space.
In particular, the hedgehog fails the following lemma
(taking $U$ to be $X$ and letting $s$ be the central point):

\begin{lemma} 
\label{lemma-fin-com}
Assume that $X$ is compact and  locally connected, 
and  $X$ has no nontrivial convergent sequences.
Fix an open $U$ with $\partial U$ finite, and
fix a finite $s \subseteq U$.
Then $U\backslash s$ has finitely many components.
\end{lemma}
\begin{proof}
Assume that $V_{n}$, for $n < \omega$, are different components
of $U\backslash s$.
Choose $x_{n} \in V_{n}$.  
Then the limit points of $\{x_{n}: n \in \omega\}$
must lie in $\partial(U\backslash s) \subseteq \partial U \cup s$. 
Thus, $\{x_{n}: n \in \omega\}$ has finitely many limit points,
which is impossible if $X$ has no nontrivial convergent sequences.
\end{proof}

We now look more closely at the locally separating points;
that is, the points $p\in X$ such that 
$U \backslash \{p\}$ is not connected
for some open connected $U \ni p$.

\begin{definition}
If $p \in U \subseteq X$, then $c(p,U)$ is the number
of components of $U \backslash \{p\}$.
\end{definition}

\begin{lemma}
\label{lemma-basic-c}
Assume that $X$ is compact and locally connected.
Fix $p \in X$ and open connected $U,V$ with $p \in U \subseteq V$.  Then:
\begin{itemizz}
\item[1.] Every component of 
$V \backslash \{p\}$ is a subset of exactly one
component of $U \backslash \{p\}$.
\item[2.] If $V \subseteq U$, then $c(p,V) \ge c(p,U)$.
\item[3.] If  $\acl(X)$ is a local base at $p$
and  $X$ has no nontrivial convergent sequences,
then $c(p,U)$ is finite.
\end{itemizz}
\end{lemma}
\begin{proof}
(1) is immediate from the fact that
if $W$ is a component of $V \backslash \{p\}$ then $W$ is
connected and $W \subseteq U \backslash \{p\}$.
For (2), use the fact that every component of $ U \backslash \{p\}$
must meet $V$ because $U$ is connected, so that (1) provides a map
from the components of  $V \backslash \{p\}$ \emph{onto}
the components of  $U \backslash \{p\}$.
For (3), choose $V \in \acl(X)$ with $p \in V \subseteq U$,
and apply (2) and Lemma \ref{lemma-fin-com}.
\end{proof}

The next lemma is trivial, but useful when $\partial U$ is finite.

\begin{lemma}
\label{lemma-all-or-none}
Suppose that $E \subseteq X$ is connected, $U \subseteq X$ is open,
and $\partial U \cap E = \emptyset$.  Then
$E \subseteq U$ or $E \cap U = \emptyset$.
\end{lemma}
\begin{proof}
$U \cap E = \overline U \cap E$ is relatively clopen in $E$,
so $U \cap E$ is either $E$ or $\emptyset$.
\end{proof}

\begin{lemma}
\label{lemma-finite-c}
Assume that $X$ is compact and locally connected,
$\acl(X)$ is a local base at $p \in X$,
and  $X$ has no nontrivial convergent sequences.  Then 
there is an $n \in \omega$ such that $ c(p,U) \le n$ for all
open connected $U \ni p$.
\end{lemma}
\begin{proof}
If this fails, then applying Lemma \ref{lemma-basic-c},
we may fix open connected $U_n \ni p$ for $n \in \omega$ such that
$U_0 \supseteq \overline U_1 \supseteq U_1 \supseteq \overline U_2 \cdots$ and
$2 \le c(p, U_0) < c(p, U_1) < \cdots$.
Then, we may define a subtree $T \subseteq \omega^{< \omega}$ 
and open connected $W_s$ for $s \in T$
and $k_s \in \omega\backslash\{0\}$ for $s \in T$ as follows:
\begin{itemizz}
\item[1.] $W_{()}$  is the component of $p$ in $X$.
\item[2.] If $\lh(s) = n$, then $k_s$ is the number of components
of $U_n \backslash \{p\}$ which are subsets of $W_s$, and these
components are listed as $\{W_{s\cat i} : i < k_s\}$.
\item[3.] $s \cat i \in T$ iff $s \in T$ and $i < k_s$.
\end{itemizz}
Item (1) is a bit artificial, but it gives $T$ a root node $()$.
For the levels below the root,
note that $| T \cap \omega^{n+1}| = c(p,U_n)$, and the
$W_s$ for $s \in T \cap \omega^{n+1}$ list the components of
$U_n \backslash \{p\}$.
Let $P(T) = \{f \in \omega^\omega: \forall n\,[ f\res n \in T]\}$ be the set of 
paths through $T$.  Since every node in $T$ has at least one child,
$|P(T)|$ is either $\aleph_0$ or $2^{\aleph_0}$.
Note that $\cl(W_{s\cat i}) \subseteq W_s \cup \{p\}$, since if
$n = \lh(s) > 0$ and $q \in \cl(W_{s\cat i})  \backslash \{p\}$, then
$q$ and the points of $W_{s\cat i}$ must all lie in the same
component of $U_{n-1} \backslash \{p\}$, which is $W_s$.

Let $H = \bigcap_n U_n = \bigcap_n \overline{U_n}$.  Then $H$ is a connected
closed $G_\delta$ containing $p$, and $H$ must be infinite, since
$p$ must have uncountable character.
For each $f \in P(T)$, let
$K_f = \bigcap_n \cl(W_{f\res n}) = \{p\}  \cup \bigcap_n W_{f\res n}$.
Then the $K_f$ are connected and infinite, since $\{p\}$
cannot be a decreasing intersection of $\omega$ infinite closed sets
(or there would be a convergent sequence).
Observe that $K_f \cap K_g = \{p\}$ whenever $f \ne g$.
Thus, if $p \in V \in \acl(X)$ then $K_f \subseteq V$ for
all but finitely many $f \in P(T)$, since $K_f \subseteq V$
whenever $K_f \cap \partial V = \emptyset$ by Lemma \ref{lemma-all-or-none}.
Now let $f_i$, for $i \in \omega$ be distinct elements of $P(T)$,
and choose $q_i \in K_{f_i} \backslash \{p\}$.
Then every neighborhood of $p$ contains all but finitely
many $q_i$, so the $q_i$ converge to $p$, a contradiction.
\end{proof}

\begin{definition}
\label{def-c}
Assume that $X$ is compact and  locally connected,  $\acl(X)$ is a base for $X$,
and  $X$ has no nontrivial convergent sequences.  Then for each $p\in X$,
define $c(p) \in \omega$ to be the largest $c(p,U)$
among all open connected $U \ni p$.
\end{definition}

By a standard chaining argument:

\begin{lemma}
\label{lemma-chain}
Assume that $X$ is compact and  locally connected and
 $\acl(X)$ is a base for $X$.  Fix a connected
open $U \subseteq X$ and a compact $F \subseteq U$.
Then there is a connected $V \in \acl(X)$ such that 
$F \subseteq V \subseteq \overline V \subseteq U$.
\end{lemma}
\begin{proof}
Let $\GG = \{W \in \acl(X): \emptyset \ne \overline W \subseteq U \ \&\ 
W \mbox{ is connected}\}$.  Then $\bigcup \GG = U$.
View $\GG$ as an undirected graph, by putting an edge between
$W_1$ and $W_2$ iff $W_1 \cap W_2 \ne \emptyset$.
Then $\GG$ is connected as a graph because
$U$ is connected and the components of $\GG$
yield topological components of $U$.
Fix a finite $\GG_0 \subseteq \GG$ such that $F \subseteq \bigcup \GG_0$.
Then fix a finite connected $\GG_1$ with
$\GG_0 \subseteq \GG_1 \subseteq \GG$. 
Let $V = \bigvee \GG_1 = \intr(\cl(\bigcup \GG_1))$.
\end{proof}

\begin{lemma}
\label{lemma-decreasing-seq}
Assume that $X$ is compact and locally connected,
$\acl(X)$ is a base for $X$,
and  $X$ has no nontrivial convergent sequences.
Then there is
no sequence of open sets $\langle U_n : n \in \omega \rangle$
such that $\overline U_{n+1} \subsetneqq U_n$ for all $n$ and
$\overline U_n \backslash U_{n+1}$ is connected for all \textnormal{even} $n$.
\end{lemma}
\begin{proof}
Given such a sequence,
choose $x_n \in \overline U_n \backslash U_{n+1}$, and
let $y$ be a limit point of $\{x_{2m} : m \in \omega \}$.
Since $\langle  x_{2m} : m \in \omega \rangle$ cannot converge to $y$,
fix a connected
$W \in \acl(X)$ and disjoint infinite $A,B \subseteq \{2m : m \in \omega\}$
such that $x_n \in W$ for all $n \in A$ and $x_n \notin W$ for all $n \in B$.
Since $\partial W$ is finite, we may also assume 
(shrinking $A,B$ if necessary) that 
$\partial W \cap (\overline U_n \backslash U_{n+1}) = \emptyset$
for all $n \in A \cup B$.
Then, by Lemma \ref{lemma-all-or-none}, 
$ \overline U_n \backslash U_{n+1} \subseteq W$ for all $n \in A$
and $ (\overline U_n \backslash U_{n+1}) \cap W = \emptyset$ for all $n \in B$.
But then, for $n \in B$,
the connected $W$ is partitioned into the disjoint open sets
$W \cap U_{n+1}$, $W \setminus \overline{U_{n}}$,
both of which are non-empty when $n > \min(A)$.
\end{proof}

\begin{lemma}
\label{lemma-allLS}
Assume that $X$ is compact and locally connected,
$\acl(X)$ is a base for $X$,
and  $X$ has no nontrivial convergent sequences.
Then every non-isolated point in $X$ is locally separating.
\end{lemma}
\begin{proof}
Suppose we have a non-isolated $p$ which is not locally separating;
so $U \backslash \{p\}$ is connected whenever $U$ is open and connected.
Then inductively construct $U_n$ for $n \in \omega$ such that
\begin{itemizz} 
\item[1.] Each $U_n$ is open and $p \in U_n$. \samepage %% don't break this
\item[2.] Each $\overline U_{n+1} \subsetneqq U_n$.\samepage %% don't break this
\item[3.] $\overline U_n \backslash U_{n+1}$ is connected 
whenever $n$ is even.  \samepage %% don't break this
\item[4.] Each $U_n \in \acl(X)$.  \samepage %% don't break this
\item[5.] $U_n$ is connected for all even $n$.
\end{itemizz}
Then (1)(2)(3) contradict Lemma \ref{lemma-decreasing-seq}.

To construct the $U_n$:  Let $U_0 \in \acl(X)$ be such that
$p \in U_0$ and $U_0$ is connected and not clopen.
Given $U_n$, where $n$ is even,
we construct $U_{n+1}$ and $U_{n+2}$ as follows:

Say $\partial U_n = \{q^j : j < r\}$; of course, $r$ and the $q^j$ depend
on $n$.  For each $j$, choose $V^j \in \acl(X)$ be such that
$q^j \in V^j$, $p \notin \cl(V^j)$, and $V^j$ is connected.
Also make sure that the $\overline{V^j}$ are disjoint;
then $\overline{V^j} \cap \partial U_n = \{q^j\}$.
Let $\{W^j_i : i < c^j\}$ list the components of $V^j \backslash \{q^j\}$;
so $2 \le c^j < \omega$.
Then $W^j_i$ is connected and $\partial U_n \cap W^j_i = \emptyset$, so
$W^j_i \subseteq U_n $ or $W^j_i \cap U_n =\emptyset$;
say $W^j_i \subseteq U_n $ for $i < d^j$ and
$W^j_i \cap U_n =\emptyset$ for $d^j \le i < c^j$; so $1 \le d^j < c^j$.
Choose $y^j_i \in W^j_i$.  Now  $U_n$ is connected and
$p$ is not locally separating, so $U_n \backslash \{p\}$ is connected.
Applying Lemma \ref{lemma-chain}, fix a connected $R \in \acl(X)$
such that $\{y^j_i : j < r \ \&\ i < d^j\} \subseteq R \subseteq \overline R
\subseteq U \backslash \{p\}$.
Let $S$ be the finite union $R \cup \bigcup \{W^j_i : j < r \ \&\ i < d^j\}$.
Then $S$ is open and connected,  $p \notin \overline S$,
and each $q^j \in \overline S$.
Let $U_{n+1} = U_n \backslash \overline S = 
\overline{U_n} \setminus \overline S$.
Then $p \in U_{n+1} \in \acl(X)$,
and $\overline U_n \backslash U_{n+1} = \overline S$ is connected.
Also, each $q^j \notin \overline{U_{n + 1}}$ 
because $ U_{n+1} \cap V^j = \emptyset$,
so that $\overline{U_{n + 1}} \subseteq U_n$.

Now, choose a connected $U_{n+2} \in \acl(X)$ so that 
$p \in U_{n+2} \subseteq \overline{U_{n+2}} \subsetneqq U_{n+1}$.
\end{proof}

\begin{theorem}
\label{thm-acl-not-base}
If $X$ is infinite, compact, locally connected,
and $\acl(X)$ is a base for $X$,
then $X$ has a nontrivial convergent sequence.
\end{theorem}
\begin{proof}
Suppose not.
Fix any non-isolated $p \in X$; then $p$
is locally separating by Lemma \ref{lemma-allLS},
so $c(p) \ge 2$ (see Definition \ref{def-c}).
Fix a connected $U \in \acl(X)$ such that $p \in U$
and $c(p,U) = c(p)$.  Let $W_i$, for $i < c(p)$ be the components
of $U \backslash \{p\}$.  Then $c(p,V) = c(p)$ whenever
$V \in \acl(X)$ and $p \in V \subseteq U$;
furthermore, the components of $V \backslash \{p\}$ are the
sets $W_i \cap V$ for $i < c(p)$.

Let $Y = \cl(W_0)$.  Then $\acl(Y)$ is a base for $Y$,
$Y$ is locally connected,
and  $Y$ has no nontrivial convergent sequences.
Furthermore, $p \in Y$ and $p$ is not locally separating in $Y$,
contradicting Lemma \ref{lemma-allLS} applied to $Y$.
\end{proof}
\section{Further Remarks}
\label{sec-rem}

We note that in constructing a locally 
connected compactum, the monotone bonding maps, as used also
by van Mill \cite{vM}, are inevitable:

\begin{remark}
Assume that $X \subseteq [0,1]^{\omega_1}$ is compact
and locally connected.
Define $X_\alpha = \pi^{\omega_1}_\alpha(X) \subseteq [0,1]^\alpha$.
Then there is a club $C \subseteq \omega_1$
such that $X_\alpha$ is locally connected for all $\alpha \in C$,
and such that
$\sigma^\beta_\alpha := \pi^\beta_\alpha\res X_\beta$
is monotone whenever $\alpha < \beta$
and $\alpha,\beta \in C \cup \{\omega_1\}$.
\end{remark}
\begin{proof}
Let $\BB$ be the family of all connected open $F_\sigma$ subsets 
of $X$.  Then $\BB$ is a base for $X$.
For $\alpha < \omega_1$, let $\BB_\alpha$ be the family of
all open $U \subseteq X_\alpha$ such that
$(\sigma^{\omega_1}_\alpha)\iv(U) \in \BB$.
Observe that each $U \in \BB_\alpha$ is connected.
Put $\alpha \in C$ iff $\BB_\alpha$ is a base for $X_\alpha$.
Then $C$ is club.

Now, it is sufficient to show that $(\sigma^{\omega_1}_\alpha)\iv\{x\}$
is connected whenever $\alpha \in C$ and $x \in X_\alpha$.
Choose $U_n \in \BB_\alpha$ with $x \in U_n \supseteq \overline{U_{n+1}}$
for all $n \in \omega$ and $\{x\} = \bigcap_n U_n = \bigcap_n \overline{U_n}$.
Each $(\sigma^{\omega_1}_\alpha)\iv(U_n)$ is in  $\BB$, so it and
its closure are connected, 
and $\cl( (\sigma^{\omega_1}_\alpha)\iv(U_{n+1}) ) \subseteq
(\sigma^{\omega_1}_\alpha)\iv(\overline{U_{n+1}}) \subseteq  
(\sigma^{\omega_1}_\alpha)\iv(U_n) $, so that
$(\sigma^{\omega_1}_\alpha)\iv\{x\}$ is the decreasing intersection
of the connected closed sets  $\cl( (\sigma^{\omega_1}_\alpha)\iv(U_{n}) )$,
and is hence connected.
\end{proof}

We do not know if conditions (C1 -- C10) in Section \ref{sec-inv} determine
$\ind(Z)$, but a minor addition to the 
construction will ensure that $Z$ 
does not have small
\emph{transfinite inductive dimension}; that is, 
$\trind(Z)=\infty$ (and hence $\ind(Z) = \infty$). 
The transfinite inductive dimension
$\trind$ is the natural generalization
of $\ind$; see \cite{ENGd}. 

\begin{theorem}
\label{thm-main+}
Assuming $\diamondsuit$, there is a locally connected
HS continuum $Z$ such that
$\dim(Z) = 1$, $\trind(Z)=\infty$,
 and $Z$ has no nontrivial convergent sequences.
\end{theorem}
To do this, we make sure that all perfect subsets are $G_\delta$ sets.
Observe that by local connectedness, every non-empty closed $G_\delta$
contains a non-empty connected closed $G_\delta$ subset, which in our $Z$
cannot be a singleton.
So, no non-empty closed $G_\delta$ can have dimension $0$.

\begin{lemma}
Assume that $X$ is compact, connected, and infinite, 
and all perfect subsets of $X$ are $G_\delta$ sets.
Assume also
that $\cchi(x,X) > \aleph_0$ for all $x \in X$, and that
in $X$, every non-empty closed $G_\delta$ set contains a
non-empty closed connected $G_\delta$ subset.  Then $\trind(X) = \infty$.
\end{lemma}

\begin{proof}
We prove by induction on ordinals $\alpha$
that $\neg [\trind(X) \le \alpha]$ for all such $X$.
This is obvious for $\alpha = 0$.
Assume $\alpha > 0$ and the inductive hypothesis holds
for all ordinals $\xi < \alpha$.
Suppose that $\trind(X) \le \alpha$.  Then there is a regular
open set $U$ such that $U \ne \emptyset$, $U \ne X$,
and $\trind (\partial U) = \xi < \alpha$.
Let $V = X \backslash \overline U$; then
$\overline U$ and $\overline V$ are perfect, so
$\partial U = \overline U \cap \overline V$ is a $G_\delta$,
and hence contains a non-empty closed connected $G_\delta$ subset $Y$.
Then $\trind(Y) \le \trind(\partial U) \le \xi$.
Since $Y$ satisfies the conditions of the lemma,
this is a contradiction.
\end{proof}

By the same argument, this space is
\emph{weird} in the sense of \cite{HK2};
that is, no perfect subset is totally disconnected.

To construct our $Z$ so that
perfect sets are $G_\delta$, we observe first
that if $Q \subseteq \MS \times [0,1]^{\omega_1}$ is perfect,
then $C := \{\alpha < \omega_1 : \pi^{\omega_1}_\alpha(Q)
\mbox{ is perfect}\}$ is a club.
One might then use $\diamondsuit$, as in \cite{DK}, to capture perfect subsets
of $Z$, but this is not necessary, since we already know that
$Z$ is HS, and we are already capturing countable sequences.
Thus, we get:

\begin{conditions}
We have $P_\alpha$ and $\PP_\alpha$ for $\alpha < \omega_1$ such that:
\begin{itemizz}
\item[C11.] $P_\alpha = \cl ( Z_\alpha \cap \{q^n_\alpha : n \in \omega\} )$
whenever  $\alpha \ge \omega$ and this set is perfect;
otherwise, $P_\alpha = Z_\alpha$. \samepage %% don't break this
\item[C12.] $\PP_\alpha =
\{ (\sigma^\alpha_\delta)\iv(P_\delta) : \delta \le \alpha \}$.
\item[C13.] 
$\sigma^{\alpha+1}_\alpha \rest ((\sigma^{\alpha+1}_\alpha)\iv(P))
\ :\  (\sigma^{\alpha+1}_\alpha)\iv(P) \onto P$ is irreducible 
for each $P \in \PP_\alpha$.
\end{itemizz}
\end{conditions}

\begin{proofof}{Theorem \ref{thm-main+}}
To obtain these conditions, note that (C13) is trivial for $P$
unless $t_\alpha \in P$.
If $t_\alpha \in P$, then, since $P$ is perfect, we may
choose a sequence of distinct points $\langle p_n : n \in \omega \rangle$
from $P \backslash \{t_\alpha\}$ converging to $t_\alpha$.  
Then, while we are accomplishing (C8), we make sure that
all points of $(\sigma^{\alpha+1}_\alpha)\iv\{t_\alpha\}$ are
(strong) limit points of the set of singletons,
$\{(\sigma^{\alpha+1}_\alpha)\iv\{p_n\} : n \in \omega \}$;
this implies irreducibility.

Now, we prove by induction on $\beta \ge \alpha$ that
$\sigma^{\beta}_\alpha \rest ((\sigma^{\beta}_\alpha)\iv(P))
\ :\  (\sigma^{\beta}_\alpha)\iv(P) \onto P$ is irreducible 
for each $P \in \PP_\alpha$.
Then, if $Q \subseteq Z$ is perfect, we use HS and (C7) to fix some 
$\alpha < \omega_1$ such that $P_\alpha = \sigma^{\omega_1}_\alpha(Q)$
and $P_\alpha$ is perfect.
Irreducibility then implies that
$Q = (\sigma^{\omega_1}_\alpha)\iv(P_\alpha)$, which is a $G_\delta$.
\end{proofof}

Finally, we remark that
our space $Z$ is \emph{dissipated} in the sense of \cite{KU2},
since in the inverse limit, only one point $t_\alpha$ gets expanded
in passing from $Z_\alpha$ to $Z_{\alpha+1}$; the inverse projection
of every other point is a singleton.
As pointed out in \cite{KU2}, this is also true of the original
Fedorchuk S-space \cite{FED-S}, where one point $t_\alpha$ got expanded
to a pair of points; here, and in
\cite{HK2} and van Mill \cite{vM},
$t_\alpha$ gets expanded to an interval.

\end{document}